\documentclass[10pt,leqno,twoside]{article}

\usepackage[T1]{fontenc}
\usepackage{lmodern}
\IfFileExists{microtype.sty}{\usepackage{microtype}}{}
\usepackage[paper=letterpaper,textwidth=5.15in,textheight=8.45in,hmarginratio=1:1,vmarginratio=1:1,headheight=15pt,headsep=16pt]{geometry}
\usepackage{amsmath,amssymb,amsthm}
\numberwithin{equation}{section}
\usepackage[hidelinks]{hyperref}
\hypersetup{
  pdftitle={The Quadratic Convex Hull of Real Stiefel Two-Frames in Dimension Three},
  pdfauthor={},
  pdfsubject={Stiefel manifold, semidefinite convex hull, matrix Fejer-Riesz theorem}
}

\makeatletter
\newcommand{\ps@journal}{%
  \let\@oddfoot\@empty
  \let\@evenfoot\@empty
  \def\@oddhead{\footnotesize\hfil\textsc{Quadratic Convex Hull of 
  St(3,2)}\hfil\llap{\thepage}}%
  \def\@evenhead{\footnotesize\rlap{\thepage}\hfil\textsc{Quadratic Convex Hull of 
  St(3,2)}\hfil}%
}
\def\@seccntformat#1{\csname the#1\endcsname.\quad}
\renewcommand{\section}{\@startsection{section}{1}{\z@}%
  {-1.85em \@plus -.4em \@minus -.2em}{.65em}{\normalfont\normalsize\bfseries}}
\renewcommand{\subsection}{\@startsection{subsection}{2}{\z@}%
  {-1.3em \@plus -.3em \@minus -.15em}{.45em}{\normalfont\normalsize\bfseries}}
\makeatother
\newtheoremstyle{journalthm}
  {9pt}{9pt}{\itshape}{}{\scshape}{.}{0.55em}
  {\thmname{#1}\ \thmnumber{#2}\thmnote{\ {\normalfont(#3)}}}
\theoremstyle{journalthm}
\newtheorem{theorem}{Theorem}[section]
\newtheorem{lemma}[theorem]{Lemma}

\newtheoremstyle{journalthmstar}
  {9pt}{9pt}{\itshape}{}{\scshape}{.}{0.55em}
  {\thmname{#1}\thmnote{\ {\normalfont(#3)}}}
\theoremstyle{journalthmstar}
\newtheorem*{theorem*}{Theorem}

\makeatletter
\renewenvironment{proof}[1][\proofname]{\par
  \pushQED{\qed}%
  \normalfont \topsep6\p@\@plus6\p@\relax
  \trivlist
  \item[\hskip\labelsep\itshape #1\@addpunct{.}]\ignorespaces
}{%
  \popQED\endtrivlist\@endpefalse
}
\makeatother

\allowdisplaybreaks[2]

\renewcommand{\le}{\leqslant}
\renewcommand{\ge}{\geqslant}
\newcommand{\R}{\mathbb{R}}
\newcommand{\C}{\mathbb{C}}
\newcommand{\Sset}{\mathcal{S}}
\newcommand{\T}{^\top}
\newcommand{\E}{\mathbb{E}}
\newcommand{\norm}[1]{\|#1\|}
\newcommand{\ip}[2]{\langle#1,#2\rangle}
\DeclareMathOperator{\tr}{tr}
\DeclareMathOperator{\rank}{rank}
\DeclareMathOperator{\conv}{conv}
\DeclareMathOperator{\cone}{cone}
\DeclareMathOperator{\aff}{aff}
\DeclareMathOperator{\adj}{adj}
\DeclareMathOperator{\im}{im}
\newcommand{\Sym}{\mathbb{S}}
\DeclareMathOperator{\diag}{diag}
\DeclareMathOperator{\St}{St}

\begin{document}
\thispagestyle{plain}

\begin{center}
  {\large\bfseries The Quadratic Convex Hull of 
  Real \(3\times2\) Stiefel Matrices\par}
\end{center}
\vspace{8pt}

\centerline{\large Zijian Zhou, Hanzhi Chen, Yong Xia}
\medskip
\centerline{School of Mathematical Sciences, Beihang University, Beijing, China}
\vspace{12pt}

\begin{quote}\small
\noindent\textbf{Abstract.}\hspace{0.45em}We prove the semidefinite convex-hull identity conjectured by
Dodig, Sto\v{s}i\'c, and Xavier
(Linear Algebra Appl.\ 475 (2015), 251--264)
for real \(3\times2\) Stiefel matrices.
The result gives an exact semidefinite reformulation of every
homogeneous quadratic optimization problem over this manifold.
The proof reduces extreme points of the proposed spectrahedron
to a circle model, whose convex hull is characterized using
the matrix Fej\'er--Riesz theorem.

\medskip
\noindent\textbf{Key words.}\hspace{0.45em}Stiefel manifold, quadratic convex hull, semidefinite programming, matrix Fej\'er--Riesz theorem, extreme points

\smallskip
\noindent\textbf{MSC codes.}\hspace{0.45em}90C22, 52A20, 15A45, 42A05
\end{quote}
\vspace{6pt}

\section{Introduction}

\subsection{Background and motivation}

The real Stiefel manifold consists of matrices with orthonormal
columns. Quadratic optimization over this manifold has a simple
formulation but a nonconvex feasible set. A standard approach is
to lift the matrix variable to a rank-one positive semidefinite
matrix, so that the quadratic objective becomes linear.
The resulting optimization problem is then determined by the
convex hull of the lifted feasible set. An exact semidefinite
description of this convex hull provides a convex reformulation
valid for every homogeneous quadratic objective.

This paper concerns real \(3\times2\) Stiefel matrices. We use the
notation
\[
\St(3,2)
=
\{Q\in\R^{3\times2}:Q\T Q=I_2\}.
\]
Writing \(Q=[a,b]\) and \(q=(a;b)\in\R^6\), we consider
\begin{equation}\label{eq:quadratic-problem}
\min_{Q=[a,b]\in\St(3,2)} q\T Cq,
\quad C\in\Sym^6,
\end{equation}
where \((a;b)\) denotes vertical concatenation and \(\Sym^6\)
is the space of real symmetric \(6\times6\) matrices.
The lifting \(X=qq\T\) gives
\(q\T Cq=\tr(CX)\).
Consequently, an exact description of the quadratic convex hull
\[
\conv\Sset,
\quad
\Sset
=
\{qq\T:q=(a;b),\ [a,b]\in\St(3,2)\},
\]
yields an equivalent convex formulation of
\eqref{eq:quadratic-problem}.

Dodig, Sto\v{s}i\'c, and Xavier \cite{DSX} proposed a semidefinite
description of \(\conv\Sset\).
In addition to the lifted orthogonality constraints, their
formulation imposes a \(4\times4\) positive semidefinite constraint
derived from the unit normal \(a\times b\) to the frame.
They proved that these conditions are necessary and conjectured
that they are sufficient. They also established sufficiency
when the vector \(v(X)\) defined below has norm one and when
\(X_{11}+X_{22}\) has rank two.
For general quadratic optimization over the Stiefel manifold,
Burer and Park \cite{BP} developed a strengthened semidefinite
relaxation based on Kronecker-product constraints.
The question addressed here is whether the formulation proposed
in \cite{DSX} is exact for every homogeneous quadratic objective
over \(\St(3,2)\).

\subsection{Main result}

For \(X\in\Sym^6\), write
\[
X=
\begin{pmatrix}
X_{11}&X_{12}\\
X_{12}\T&X_{22}
\end{pmatrix},
\quad X_{ij}\in\R^{3\times3}.
\]
Define the affine subspace
\begin{equation}\label{eq:affineA}
\mathcal A
=
\left\{
X\in\Sym^6:
\tr X_{11}
=
\tr X_{22}=1,\ 
\tr X_{12}=0
\right\}.
\end{equation}
These trace equations express the orthonormality of the two
frame vectors after lifting.

Set \(D=X_{11}+X_{22}\) and define
\[
v=v(X)
=
\begin{pmatrix}
(X_{12})_{23}-(X_{12})_{32}\\
(X_{12})_{31}-(X_{12})_{13}\\
(X_{12})_{12}-(X_{12})_{21}
\end{pmatrix}.
\]
The candidate spectrahedron is
\begin{equation}\label{eq:sets}
\Sigma
=
\{X\in\mathcal A:X\succeq0,\ L(X)\succeq0\},
\quad
L(X)
=
\begin{pmatrix}
1&v\T\\
v&I_3-D
\end{pmatrix}.
\end{equation}
This is precisely the formulation in
\cite[(11)--(15)]{DSX}, up to a simultaneous permutation of the
rows and columns of the \(4\times4\) block.

The additional semidefinite constraint has a direct geometric
meaning. For a lifted frame \(X=qq\T\), let \(n=a\times b\).
Then \(v(X)=n\) and \(I_3-D=nn\T\), so
\[
L(qq\T)
=
\begin{pmatrix}1\\n\end{pmatrix}
\begin{pmatrix}1\\n\end{pmatrix}^{\!\top}
\succeq0.
\]
Thus every lifted frame belongs to \(\Sigma\), and hence
\(\conv\Sset\subseteq\Sigma\).
Our main result establishes the reverse inclusion.

\begin{theorem}[Resolution of the
Dodig--Sto\v{s}i\'c--Xavier conjecture]\label{thm:main}
The quadratic convex hull of real \(3\times2\) Stiefel matrices
satisfies
\[
\conv\Sset=\Sigma.
\]
\end{theorem}

In particular, for every \(C\in\Sym^6\),
\[
\min_{Q=[a,b]\in\St(3,2)}
(a;b)\T C(a;b)
=
\min_{X\in\Sigma}\tr(CX).
\]
The semidefinite formulation proposed in \cite{DSX} is therefore
exact for every homogeneous quadratic objective.
Moreover, every matrix in \(\Sigma\) is a finite convex combination
of lifted Stiefel frames.

\subsection{Proof outline}

The proof separates the extreme-point argument from the
decomposition of matrices with singular \(L(X)\).
We first show that an extreme point of \(\Sigma\) cannot satisfy
\(L(X)\succ0\): this strict inequality forces
\(\rank X\ge3\), whereas the Barvinok--Pataki rank
bound applied to the three trace constraints would give
\(\rank X\le2\).

We then show that every feasible matrix with singular \(L(X)\)
belongs to \(\conv\Sset\).
A suitable kernel vector reduces the problem to frames whose
unit normals lie on a latitude circle.
An invertible real congruence identifies the resulting affine
slice with an auxiliary circle model, whose convex hull is
characterized using conic duality and the matrix
Fej\'er--Riesz theorem.
Compactness of \(\Sigma\) then completes the proof.

Section~2 develops the circle model, proves the latitude
reduction, and completes the extreme-point argument.
Section~3 concludes the paper.

\section{Proof of Theorem~\ref{thm:main}}

We first record the frame identities underlying the
normal-moment constraint. For \(Q=[a,b]\in\St(3,2)\), set
\(q=(a;b)\) and \(n=a\times b\). Since \(X_{12}=ab\T\)
for \(X=qq\T\), the vector triple product gives
\begin{equation}\label{eq:cross}
(X_{12}\T-X_{12})x
=
(ba\T-ab\T)x
=
n\times x
\quad (x\in\R^3).
\end{equation}
Thus \(v(qq\T)=n\). Here and below, \([u]_\times\) denotes
the skew-symmetric matrix satisfying \([u]_\times x=u\times x\),
so \(X_{12}\T-X_{12}=[v(X)]_\times\) for every
\(X\in\Sym^6\).

For an arbitrary \(X\in\Sym^6\), recall that
\(D=X_{11}+X_{22}\) and \(v=v(X)\). Define the Schur complement
of the top-left entry of \(L(X)\) by
\begin{equation}\label{eq:LPhi}
\Phi(X)=I_3-D-vv\T.
\end{equation}
Since that entry is \(1\), we have
\[
L(X)\succeq0
\iff \Phi(X)\succeq0,
\quad
L(X)\succ0
\iff \Phi(X)\succ0.
\]

The block \(L(X)\) has a moment interpretation.
Consider a finite probability distribution of frames and let
\(X=\E[qq\T]\). By linearity,
\(v=\E[n]\), while the identity
\(aa\T+bb\T+nn\T=I_3\) gives
\(I_3-D=\E[nn\T]\). Consequently,
\[
L(X)
=
\E\left[
\begin{pmatrix}1\\n\end{pmatrix}
\begin{pmatrix}1\\n\end{pmatrix}^{\!\top}
\right],
\quad
\Phi(X)
=
\E[nn\T]-\E[n]\E[n]\T.
\]
In particular, \(L(X)\) is positive semidefinite and
\(\Phi(X)\) is the covariance matrix of the unit normal.

This interpretation also explains the latitude reduction
used below. For \(w\in\R^3\) and \(\nu\in\R\),
\[
\begin{pmatrix}-\nu\\w\end{pmatrix}^{\!\top}
L(X)
\begin{pmatrix}-\nu\\w\end{pmatrix}
=
\E[(w\T n-\nu)^2].
\]
Hence, if \((-\nu;w)\in\ker L(X)\), the normal of every
positive-weight atom satisfies \(w\T n=\nu\).
When \(\norm{w}=1\) and \(|\nu|<1\), this plane intersects the
unit sphere in a nondegenerate latitude circle.
This observation motivates the construction below; the
existence of a frame decomposition will be proved without
assuming such a representation of \(X\).

The forward inclusion \(\conv\Sset\subseteq\Sigma\) follows
from the frame identities above. For the reverse inclusion,
compactness reduces the problem to the extreme points of
\(\Sigma\). Lemma~\ref{lem:extreme} shows that every extreme
point has singular \(L(X)\), and Lemma~\ref{lem:latitude}
shows that every feasible matrix with singular \(L(X)\)
belongs to \(\conv\Sset\). We begin with the auxiliary circle
model needed for the latter decomposition.

\subsection{The circle model via matrix Fej\'er--Riesz}
\label{sec:circle}

We will use the following matrix version of the Fej\'er--Riesz theorem.

\begin{theorem*}[Matrix Fej\'er--Riesz {\cite[Theorem~2.1]{DR}}; see also {\cite{BW}}]
Let
\[
Q(\theta) = \sum_{k=-d}^{d} C_k e^{ik\theta}
\qquad \bigl(C_{-k}=C_k^*\in\C^{m\times m}\bigr)
\]
be a Hermitian matrix trigonometric polynomial of degree at most \(d\) such that \(Q(\theta)\succeq0\) for all \(\theta\in\R\). Then there exists a matrix polynomial \(R(z)=\sum_{k=0}^{d} R_k z^k\in\C^{m\times m}[z]\) of degree at most \(d\) such that \(Q(\theta)=R(e^{i\theta})^*R(e^{i\theta})\) for all \(\theta\in\R\).
\end{theorem*}

We first introduce an auxiliary circle model, with parameters \((c,s)\in S^1\) and \(y=(y_1,y_2)\T\in S^1\). Define
\begin{align*}
k(c,s) &= (1,c,s)\T, &
F &= \diag(-1,1,1),\\
g &= \bigl(y_1 k(c,s);\, y_2 k(c,s)\bigr)\in\R^6, &
\mathcal K &= \bigl\{gg\T : c^2+s^2=1,\; y_1^2+y_2^2=1\bigr\}\subset\Sym^6.
\end{align*}
Partition any \(W=(H_{\alpha\beta})_{\alpha,\beta=1}^2\in\Sym^6\) into \(3\times3\) blocks. Let \(\mathcal B_0\subset\Sym^6\) be the linear subspace defined by
\begin{equation}\label{eq:B0}
H_{12} = H_{12}\T, \quad
\tr(FH_{\alpha\beta}) = 0 \quad (1\le\alpha\le\beta\le2),
\end{equation}
and define the trace-normalized affine slice
\begin{equation}\label{eq:circle}
\mathcal B = \bigl\{W\in\mathcal B_0 : \tr W = 2\bigr\}.
\end{equation}
Note that on \(\mathcal B_0\), the condition \(\tr(FH_{\alpha\alpha})=0\) gives \((H_{\alpha\alpha})_{22}+(H_{\alpha\alpha})_{33}=(H_{\alpha\alpha})_{11}\), so \(\tr W = 2\bigl((H_{11})_{11}+(H_{22})_{11}\bigr)\); thus the trace normalization \(\tr W=2\) is equivalent to \((H_{11})_{11}+(H_{22})_{11}=1\).

\begin{lemma}[Circle model hull]\label{lem:circle}
The sets \(\mathcal K\) and \(\mathcal B\) satisfy
\[
\conv\mathcal K = \{W\in\mathcal B : W\succeq0\},
\]
and \(\aff\mathcal K=\mathcal B\) has dimension \(14\).
\end{lemma}

\begin{proof}
We organize the proof in four steps: Step~1 verifies the forward inclusion and uses conic duality to reduce the problem to an inclusion of dual cones; Step~2 constructs a positive semidefinite certificate via matrix Fej\'er--Riesz factorization; Step~3 characterizes the residual in \(\mathcal B_0^\perp\) to establish cone equality; and Step~4 intersects with the trace hyperplane and determines the affine hull and its dimension.

\emph{Step 1: Forward inclusion and reduction to dual cones.}
For each \(gg\T\in\mathcal K\), the \(3\times3\) blocks are \(H_{\alpha\beta}=y_\alpha y_\beta k(c,s)k(c,s)\T\). Hence \(H_{12}=H_{12}\T\), \(\tr(FH_{\alpha\beta})=y_\alpha y_\beta\bigl(k(c,s)\T Fk(c,s)\bigr)=0\) since \(c^2+s^2=1\), and
\begin{equation}\label{eq:circle-atom-trace}
\tr(gg\T)=\norm{g}^2=(y_1^2+y_2^2)(1+c^2+s^2)=2.
\end{equation}
Thus \(\conv\mathcal K\subseteq\{W\in\mathcal B:W\succeq0\}\).

Introduce the convex cones
\begin{align*}
\mathcal C_{\mathcal K} &= \cone(\mathcal K) = \{tW : t\ge0,\ W\in\conv\mathcal K\},\\
\mathcal C_{\mathcal B} &= \{W\in\mathcal B_0 : W\succeq0\}.
\end{align*}
From the forward inclusion we have \(\mathcal C_{\mathcal K}\subseteq\mathcal C_{\mathcal B}\), and \(\mathcal C_{\mathcal B}\) is closed. Because \(\mathcal K\) is the continuous image of the compact set \(S^1\times S^1\), it is compact, so \(\conv\mathcal K\) is compact \cite[Theorem~17.2]{Rockafellar}. By \eqref{eq:circle-atom-trace}, every matrix in \(\conv\mathcal K\) has trace \(2\), so \(0\notin\conv\mathcal K\), whence \(\mathcal C_{\mathcal K}\) is also closed \cite[Corollary~9.6.1]{Rockafellar}.

Equip \(\Sym^6\) with the Frobenius inner product \(\ip{A}{W}=\tr(AW)\), and define the dual cone of any \(\mathcal C\subseteq\Sym^6\) by
\[
\mathcal C^* = \{A\in\Sym^6 : \ip{A}{W}\ge0 \text{ for all } W\in\mathcal C\}.
\]
Since \(\mathcal C_{\mathcal K}\subseteq\mathcal C_{\mathcal B}\), taking dual cones reverses the inclusion to \(\mathcal C_{\mathcal B}^*\subseteq\mathcal C_{\mathcal K}^*\). Because both \(\mathcal C_{\mathcal K}\) and \(\mathcal C_{\mathcal B}\) are closed convex cones, by the bipolar theorem \cite[Theorem~14.1]{Rockafellar} it therefore suffices to prove \(\mathcal C_{\mathcal K}^*\subseteq\mathcal C_{\mathcal B}^*\).

\emph{Step 2: Sum-of-squares certificate via matrix Fej\'er--Riesz.}
Let \(A=(A_{\alpha\beta})_{\alpha,\beta=1}^2\in\mathcal C_{\mathcal K}^*\). We will construct \(P\succeq0\) such that \(A-P\in\mathcal B_0^\perp\).

For any \(\theta\in\R\), set \(k_\theta:=k(\cos\theta,\sin\theta)=(1,\cos\theta,\sin\theta)\T\) and define \(Q_A(\theta)\in\Sym^2\) entrywise by
\begin{equation}\label{eq:QA-def}
(Q_A(\theta))_{\alpha\beta} := k_\theta\T \widetilde A_{\alpha\beta} k_\theta,
\qquad
\widetilde A_{\alpha\beta} := \tfrac12(A_{\alpha\beta}+A_{\alpha\beta}\T)
\quad (\alpha,\beta\in\{1,2\}).
\end{equation}
Since \(A\in\Sym^6\) satisfies \(A_{21}=A_{12}\T\), we have \(\widetilde A_{21}=\widetilde A_{12}\), so \(Q_A(\theta)\) is symmetric. Now take any \(x=(x_1,x_2)\T\in\R^2\) and put \(g=(x_1 k_\theta;\,x_2 k_\theta)\). If \(x=0\), then \(gg\T=0\in\mathcal C_{\mathcal K}\). If \(x\ne0\), then \(x/\norm{x}\in S^1\), so \(\hat g:=g/\norm{x}\) yields a generator \(\hat g\hat g\T\in\mathcal K\), and hence \(gg\T=\norm{x}^2\hat g\hat g\T\in\mathcal C_{\mathcal K}\). Because \(A\in\mathcal C_{\mathcal K}^*\), it follows that
\begin{equation}\label{eq:QA-psd}
0 \le \ip{A}{gg\T} = \tr(Agg\T) = g\T A g = x\T Q_A(\theta)x.
\end{equation}
By \eqref{eq:QA-psd}, \(Q_A(\theta)\) is positive semidefinite on \(\R^2\) for every \(\theta\in\R\). For any \(z=u+iv\in\C^2\) (\(u,v\in\R^2\)), symmetry of \(Q_A(\theta)\) gives \(z^*Q_A(\theta)z=u\T Q_A(\theta)u+v\T Q_A(\theta)v\ge0\), so \(Q_A(\theta)\) is also positive semidefinite as a Hermitian matrix on \(\C^2\).

Substituting \(\cos\theta=(e^{i\theta}+e^{-i\theta})/2\) and \(\sin\theta=(e^{i\theta}-e^{-i\theta})/(2i)\) into \eqref{eq:QA-def} expresses \(Q_A(\theta)\) as a Hermitian matrix trigonometric polynomial of degree at most \(2\):
\[
Q_A(\theta) = \sum_{k=-2}^{2} C_k e^{ik\theta},
\qquad C_{-k}=C_k^*=\overline{C_k}\in\C^{2\times2} \quad (k=0,1,2),
\]
where the \((\alpha,\beta)\) entries of \(C_0,C_1,C_2\) are given explicitly by
\begin{align*}
(C_0)_{\alpha\beta} &= (\widetilde A_{\alpha\beta})_{11}
  + \tfrac12\bigl((\widetilde A_{\alpha\beta})_{22}+(\widetilde A_{\alpha\beta})_{33}\bigr),\\
(C_1)_{\alpha\beta} &= (\widetilde A_{\alpha\beta})_{12} - i(\widetilde A_{\alpha\beta})_{13},\\
(C_2)_{\alpha\beta} &= \tfrac14\bigl((\widetilde A_{\alpha\beta})_{22}
  - (\widetilde A_{\alpha\beta})_{33} - 2i(\widetilde A_{\alpha\beta})_{23}\bigr).
\end{align*}
By the matrix Fej\'er--Riesz theorem, there exists a matrix polynomial \(R(z)=R_0+R_1 z+R_2 z^2\in\C^{2\times2}[z]\) such that \(Q_A(\theta)=R(e^{i\theta})^*R(e^{i\theta})\) for all \(\theta\in\R\). Multiplication by the unit-modulus scalar \(e^{-i\theta}\) preserves \(R^*R\) and shifts the Fourier modes from \(0,1,2\) to \(-1,0,1\). Separating
\[
\widetilde R(\theta) := e^{-i\theta}R(e^{i\theta})
= R_1 + (R_0+R_2)\cos\theta + i(R_2-R_0)\sin\theta
\]
into its real and imaginary parts \(\widetilde R(\theta)=U(\theta)+iV(\theta)\) gives
\begin{align*}
U(\theta) &= \operatorname{Re}\widetilde R(\theta)
           = \operatorname{Re}(R_1) + \operatorname{Re}(R_0+R_2)\cos\theta - \operatorname{Im}(R_2-R_0)\sin\theta,\\
V(\theta) &= \operatorname{Im}\widetilde R(\theta)
           = \operatorname{Im}(R_1) + \operatorname{Im}(R_0+R_2)\cos\theta + \operatorname{Re}(R_2-R_0)\sin\theta.
\end{align*}
Since \(Q_A(\theta)=\widetilde R(\theta)^*\widetilde R(\theta)=(U\T U+V\T V)+i(U\T V-V\T U)\) is real, taking real parts yields
\[
Q_A(\theta) = U(\theta)\T U(\theta) + V(\theta)\T V(\theta)
\qquad\text{for all }\theta\in\R.
\]
Write each \(R_j=\begin{pmatrix}r_{j,1}&r_{j,2}\end{pmatrix}\in\C^{2\times2}\) with columns \(r_{j,1},r_{j,2}\in\C^2\) (\(j=0,1,2\)), and define the \(2\times3\) real column blocks
\begin{align*}
B_{1,\alpha} &= \begin{pmatrix}
  \operatorname{Re}(r_{1,\alpha}) & \operatorname{Re}(r_{0,\alpha}+r_{2,\alpha}) & -\operatorname{Im}(r_{2,\alpha}-r_{0,\alpha})
\end{pmatrix},\\
B_{2,\alpha} &= \begin{pmatrix}
  \operatorname{Im}(r_{1,\alpha}) & \operatorname{Im}(r_{0,\alpha}+r_{2,\alpha}) & \operatorname{Re}(r_{2,\alpha}-r_{0,\alpha})
\end{pmatrix}
\qquad (\alpha=1,2).
\end{align*}
Setting \(B_1=\begin{pmatrix}B_{1,1}&B_{1,2}\end{pmatrix}\in\R^{2\times6}\) and \(B_2=\begin{pmatrix}B_{2,1}&B_{2,2}\end{pmatrix}\in\R^{2\times6}\), we have \(U(\theta)x=B_1 g\) and \(V(\theta)x=B_2 g\) for \(g=(x_1 k_\theta;\,x_2 k_\theta)\). With \(P:=B_1\T B_1+B_2\T B_2\succeq0\), we obtain
\[
g\T A g = \norm{U(\theta)x}^2 + \norm{V(\theta)x}^2
= \norm{B_1 g}^2 + \norm{B_2 g}^2
= g\T(B_1\T B_1+B_2\T B_2)g
= g\T P g
\]
for all \(\theta\in\R\) and \(x\in\R^2\).

\emph{Step 3: Residual characterization and cone equality.}
Put \(M=A-P\in\Sym^6\), and define \(Q_M(\theta)\in\Sym^2\) together with the symmetrized \(3\times3\) blocks \(\widetilde M_{\alpha\beta}\) analogously to \eqref{eq:QA-def}. Then \(x\T Q_M(\theta)x=g\T Mg=0\) for all \(x\in\R^2\) and \(\theta\in\R\). Since \(Q_M(\theta)\in\Sym^2\), it follows that \(Q_M(\theta)=0\); hence each symmetrized block \(N=\widetilde M_{\alpha\beta}\) satisfies \(k_\theta\T N k_\theta=0\) for all \(\theta\in\R\). Expanding in the linearly independent Fourier basis \(\{1,\cos\theta,\sin\theta,\cos2\theta,\sin2\theta\}\) yields
\[
\Bigl(N_{11}+\tfrac{N_{22}+N_{33}}{2}\Bigr)
+ 2N_{12}\cos\theta + 2N_{13}\sin\theta
+ \tfrac{N_{22}-N_{33}}{2}\cos2\theta + N_{23}\sin2\theta = 0,
\]
forcing \(N_{12}=N_{13}=N_{23}=0\), \(N_{22}=N_{33}\), and \(N_{11}=-N_{22}\), i.e., \(\widetilde M_{\alpha\beta}=\lambda_{\alpha\beta}F\) for some \(\lambda_{\alpha\beta}\in\R\) (\(1\le\alpha\le\beta\le2\)), where \(F=\diag(-1,1,1)\). Consequently, \(M_{11}=\widetilde M_{11}=\lambda_{11}F\) and \(M_{22}=\widetilde M_{22}=\lambda_{22}F\). For the off-diagonal block, setting \(K:=\tfrac12(M_{12}-M_{12}\T)\) gives \(K\T=-K\) and
\[
M_{12} = \tfrac12(M_{12}+M_{12}\T) + \tfrac12(M_{12}-M_{12}\T)
= \widetilde M_{12}+K = \lambda_{12}F+K,
\]
with \(M_{21}=M_{12}\T=\lambda_{12}F-K\) since \(M\in\Sym^6\).

Thus the residual \(M=A-P\) belongs to the \(6\)-dimensional linear subspace \(\mathcal M\subset\Sym^6\) consisting of all matrices of the form
\[
Z = \begin{pmatrix}
  \lambda_{11}F & \lambda_{12}F+K\\
  \lambda_{12}F-K & \lambda_{22}F
\end{pmatrix}
\qquad (\lambda_{11},\lambda_{12},\lambda_{22}\in\R,\; K\T=-K).
\]
For any \(Z\in\mathcal M\) of this form and any \(W=(H_{\alpha\beta})\in\Sym^6\), using \(F\T=F\), \(K\T=-K\), and \(H_{21}=H_{12}\T\), direct block expansion gives
\begin{equation}\label{eq:MW-inner}
\ip{Z}{W}
= \lambda_{11}\tr(FH_{11}) + 2\lambda_{12}\tr(FH_{12})
  + \lambda_{22}\tr(FH_{22}) + 2\tr(K\T H_{12}).
\end{equation}
Notice that \(2\tr(K\T H_{12})=\tr\bigl(K\T(H_{12}-H_{12}\T)\bigr)\). If \(W\in\mathcal B_0\), then \eqref{eq:MW-inner} vanishes for every \(Z\in\mathcal M\). Conversely, suppose \(\ip{Z}{W}=0\) for all \(Z\in\mathcal M\). Setting \(K=0\) and varying \(\lambda_{11},\lambda_{12},\lambda_{22}\) independently gives \(\tr(FH_{\alpha\beta})=0\) for \(1\le\alpha\le\beta\le2\). Setting \(\lambda_{11}=\lambda_{12}=\lambda_{22}=0\) and choosing the skew-symmetric matrix \(K=H_{12}-H_{12}\T\) gives
\[
0 = 2\tr(K\T H_{12}) = \tr\bigl(K\T(H_{12}-H_{12}\T)\bigr)
  = \tr(K\T K) = \norm{K}_{\mathrm F}^2,
\]
forcing \(H_{12}=H_{12}\T\). Thus \(\mathcal B_0=\mathcal M^\perp\) and \(\mathcal B_0^\perp=\mathcal M\).

In particular, since \(M\in\mathcal M=\mathcal B_0^\perp\), every \(W\in\mathcal C_{\mathcal B}\) satisfies \(\ip{M}{W}=0\) and hence
\[
\ip{A}{W} = \ip{P}{W} + \ip{M}{W}
= \tr(PW) = \tr\bigl(W^{1/2}PW^{1/2}\bigr) \ge 0.
\]
This establishes \(\mathcal C_{\mathcal K}^*\subseteq\mathcal C_{\mathcal B}^*\), whence \(\mathcal C_{\mathcal K}^*=\mathcal C_{\mathcal B}^*\), and the bipolar theorem gives \(\mathcal C_{\mathcal K}=\mathcal C_{\mathcal B}\).

\emph{Step 4: Trace slice, affine hull, and dimension.}
Intersect both sides of \(\mathcal C_{\mathcal K}=\mathcal C_{\mathcal B}\) with the hyperplane \(\{W:\tr W=2\}\). On the left, each \(W\in\mathcal C_{\mathcal K}\) has the form \(tW_0\) with \(t\ge0\) and \(W_0\in\conv\mathcal K\). By \eqref{eq:circle-atom-trace} and linearity of the trace, \(\tr W_0=2\). Thus \(\tr W=2\) is equivalent to \(t=1\), so the left slice is \(\conv\mathcal K\). On the right, by the definitions of \(\mathcal C_{\mathcal B}\) and \(\mathcal B\), the slice is \(\{W\in\mathcal B:W\succeq0\}\). Equating the two slices yields
\[
\conv\mathcal K = \{W\in\mathcal B : W\succeq0\}.
\]

Finally, \(\dim\mathcal B_0=\dim\Sym^6-\dim\mathcal M=21-6=15\). Averaging the eight lifts in \(\mathcal K\) with \((c,s)\in\{(\pm1,0),(0,\pm1)\}\) and \(y\in\{(1,0)\T,(0,1)\T\}\) with equal weight \(1/8\) gives the positive definite point
\[
W_* = \tfrac12 I_2\otimes\diag\!\bigl(1,\tfrac12,\tfrac12\bigr) \in \conv\mathcal K \subset \mathcal B.
\]
Since \(W_*\in\mathcal B_0\) has trace \(2\), the trace is a nonzero linear functional on \(\mathcal B_0\). Its level set \(\mathcal B\) therefore has codimension one, giving \(\dim\mathcal B=15-1=14\). Let \(W\in\mathcal B\) be arbitrary. Since \(\mathcal B\) is affine,
\[
W_\varepsilon = W_* + \varepsilon(W-W_*) \in \mathcal B.
\]
For sufficiently small \(\varepsilon>0\), \(W_\varepsilon\succ0\) by openness of the positive definite cone and \(W_*\succ0\). Hence both \(W_*\) and \(W_\varepsilon\) belong to \(\conv\mathcal K\) by the convex-hull equality just proved. Since
\[
W = (1-\varepsilon^{-1})W_* + \varepsilon^{-1}W_\varepsilon
\]
is an affine combination of these two points, \(W\in\aff\mathcal K\). Together with \(\mathcal K\subseteq\mathcal B\), this gives \(\aff\mathcal K=\mathcal B\) and \(\dim\aff\mathcal K=14\).
\end{proof}

\subsection{Reduction to a latitude}\label{sec:latitude}

For a unit vector \(w\in\R^3\) and a scalar \(\nu\) with \(|\nu|<1\), let
\[
\mathcal A_{\nu,w} = \Bigl\{X\in\mathcal A : L(X)\begin{pmatrix}-\nu\\ w\end{pmatrix}=0\Bigr\}.
\]

\begin{lemma}[Decomposition on singular faces]\label{lem:latitude}
If \(X\in\Sigma\) and \(L(X)\) is singular, then \(X\in\conv\Sset\).
\end{lemma}

\begin{proof}
We first fix a unit vector \(w\) and a scalar \(\nu\) with \(|\nu|<1\), and show that every positive semidefinite \(X\in\mathcal A_{\nu,w}\) belongs to \(\conv\Sset\). Step~1 rotates \(w\) to \(e_3\), preserving positive semidefiniteness and transforming the kernel condition. Step~2 proves the convex-hull inclusion for the rotated matrix by reducing to the circle model of Lemma~\ref{lem:circle}. Step~3 then shows that every \(X\in\Sigma\) with singular \(L(X)\) has a kernel vector of this form.

\emph{Step 1: Reduction to \(w=e_3\) by rotation.}
Choose \(P\in\mathrm{SO}(3)\) with \(Pw=e_3\) and set \(\widetilde X=(I_2\otimes P)X(I_2\otimes P)\T\). The congruence preserves positive semidefiniteness, so it remains to check that \(\widetilde X\in\mathcal A_{\nu,e_3}\). Each \(3\times3\) block transforms as \(\widetilde X_{ij}=PX_{ij}P\T\), which satisfies \(\tr(PX_{ij}P\T)=\tr X_{ij}\); hence the three trace constraints in \eqref{eq:affineA} are preserved and \(\widetilde X\in\mathcal A\). Since \(D=X_{11}+X_{22}\), we have \(\widetilde D=PDP\T\). Because \(\det P=1\), \(P\) preserves the cross product, so \(P(v\times P\T x)=(Pv)\times x\) for all \(x\in\R^3\), i.e.,
\[
\widetilde X_{12}\T-\widetilde X_{12}
= P(X_{12}\T-X_{12})P\T
= P[v]_\times P\T
= [Pv]_\times.
\]
Hence \(\widetilde v=Pv\) and \(L(\widetilde X)=\diag(1,P)\,L(X)\,\diag(1,P)\T\). In particular,
\[
L(\widetilde X)\begin{pmatrix}-\nu\\ e_3\end{pmatrix}
= \diag(1,P)\,L(X)\begin{pmatrix}-\nu\\ w\end{pmatrix}
= 0,
\]
showing \(\widetilde X\in\mathcal A_{\nu,e_3}\).

\emph{Step 2: Reduction of the \(e_3\)-slice to the circle model.}
We now work with \(w=e_3\), writing \(X\) for the rotated matrix. Set \(\kappa=\sqrt{1-\nu^2}>0\). For \(X\in\mathcal A\), membership in \(\mathcal A_{\nu,e_3}\) is equivalent to the kernel equation
\[
L(X)\begin{pmatrix}-\nu\\ e_3\end{pmatrix}
= \begin{pmatrix}1 & v\T\\ v & I_3-D\end{pmatrix}\begin{pmatrix}-\nu\\ e_3\end{pmatrix}
= \begin{pmatrix}v_3-\nu\\ e_3-De_3-\nu v\end{pmatrix}
= 0.
\]
The first entry gives \(v_3=\nu\). Substituting this into \(1-D_{33}-\nu v_3=0\) gives \(D_{33}=\kappa^2\), so the kernel condition is equivalent to the four linear equations
\begin{equation}\label{eq:kernel}
\begin{aligned}
v_3 &= \nu, &
D_{33} &= \kappa^2,\\
D_{13}+\nu v_1 &= 0, &
D_{23}+\nu v_2 &= 0.
\end{aligned}
\end{equation}
The three trace equations \eqref{eq:affineA} act on disjoint blocks, so \(\dim\mathcal A=21-3=18\), with direction space \(\mathcal A_0=\{H\in\Sym^6:\tr H_{11}=\tr H_{22}=\tr H_{12}=0\}\). Letting \(E_{ij}\in\R^{3\times3}\) denote the \((i,j)\) standard matrix unit, we consider four perturbations in \(\mathcal A_0\), each altering exclusively one of the four left-hand sides in \eqref{eq:kernel}:
\begin{enumerate}
\item \(E_{12}-E_{21}\) in \(X_{12}\) (with \(X_{21}=X_{12}\T\)) alters \(v_3\) by \(2\);
\item \(E_{13}+E_{31}\) in \(X_{11}\) alters \(D_{13}+\nu v_1\) by \(1\);
\item \(E_{23}+E_{32}\) in \(X_{11}\) alters \(D_{23}+\nu v_2\) by \(1\);
\item \(\diag(-1,0,1)\) in \(X_{11}\) alters \(D_{33}\) by \(1\).
\end{enumerate}
Hence the linear part of the four left-hand sides in \eqref{eq:kernel}, restricted to \(\mathcal A_0\), has rank \(4\). The map of left-hand sides in \eqref{eq:kernel} is therefore onto \(\R^4\), so the prescribed level set \(\mathcal A_{\nu,e_3}\) is nonempty and has dimension \(18-4=14\).

To apply Lemma~\ref{lem:circle} to this slice, we first map the circle-model generators to frames whose normals have third coordinate \(\nu\). For fixed \(\nu\), define \(T_\nu:\R^6\to\R^6\) by
\[
T_\nu = \begin{pmatrix}
  0 & 0 & -1 & 0 & -\nu & 0\\
  0 & 1 & 0 & 0 & 0 & -\nu\\
  0 & 0 & 0 & \kappa & 0 & 0\\
  0 & -\nu & 0 & 0 & 0 & 1\\
  0 & 0 & -\nu & 0 & -1 & 0\\
  \kappa & 0 & 0 & 0 & 0 & 0
\end{pmatrix}.
\]
Reordering the rows as \((6,3,1,5,2,4)\) and the columns as \((1,4,3,5,2,6)\) gives
\[
\diag\!\left(
  \kappa I_2,\;
  \begin{pmatrix}-1 & -\nu\\ -\nu & -1\end{pmatrix},\;
  \begin{pmatrix}1 & -\nu\\ -\nu & 1\end{pmatrix}
\right).
\]
Each block has determinant \(\kappa^2=1-\nu^2>0\), so \(T_\nu\) is invertible.

Set
\begin{align*}
\mathcal G &= \bigl\{(y_1,y_1 c,y_1 s,y_2,y_2 c,y_2 s)\T : c^2+s^2=y_1^2+y_2^2=1\bigr\},\\
\mathcal F_\nu &= \bigl\{(a;b) : [a,b]\in\St(3,2),\ (a\times b)_3=\nu\bigr\}.
\end{align*}
We verify that \(T_\nu|_{\mathcal G}:\mathcal G\to\mathcal F_\nu\) is a bijection. Take any \(g\in\mathcal G\), written as \(g=(y_1 k(c,s);\,y_2 k(c,s))\), and let \(e=(-s,c,0)\T\) and \(f=(-\nu c,-\nu s,\kappa)\T\). These two vectors form an orthonormal pair with \(e\times f=(\kappa c,\kappa s,\nu)\T\); geometrically, \((c,s)\in S^1\) locates the normal along the latitude circle, while \(y=(y_1,y_2)\T\in S^1\) rotates the frame within its tangent plane. Writing \(T_\nu g=(a;b)\), direct matrix multiplication gives \(a=y_1 e+y_2 f\) and \(b=y_1 f-y_2 e\). Thus \([a,b]\in\St(3,2)\) with normal \(a\times b=(y_1^2+y_2^2)(e\times f)=(\kappa c,\kappa s,\nu)\T\), so \(T_\nu g\in\mathcal F_\nu\).

Conversely, given \((a;b)\in\mathcal F_\nu\), let \(n=a\times b\). Since \(\norm{n}=1\) and \(n_3=\nu\), setting \(c=n_1/\kappa\) and \(s=n_2/\kappa\) gives \(c^2+s^2=1\). The corresponding orthonormal vectors \(e,f\) form a positively oriented basis of \(n^\perp\). Put \(y_1=e\T a\) and \(y_2=f\T a\); then \(y_1^2+y_2^2=1\), \(a=y_1 e+y_2 f\), and \(b=n\times a=(e\times f)\times(y_1 e+y_2 f)=y_1 f-y_2 e\), so \(T_\nu g=(a;b)\) for \(g=(y_1 k(c,s);\,y_2 k(c,s))\in\mathcal G\). Injectivity of \(T_\nu|_{\mathcal G}\) follows from invertibility of \(T_\nu\).

Now take any \(gg\T\in\mathcal K\) and put \(q=T_\nu g=(a;b)\in\mathcal F_\nu\) with \(n=a\times b\). Then \(T_\nu gg\T T_\nu\T=qq\T\in\Sset\subseteq\mathcal A\), and since \(n_3=\nu\),
\[
L(qq\T)\begin{pmatrix}-\nu\\ e_3\end{pmatrix}
= (n_3-\nu)\begin{pmatrix}1\\ n\end{pmatrix}
= 0,
\]
so \(T_\nu gg\T T_\nu\T\in\mathcal A_{\nu,e_3}\). The map \(W\mapsto T_\nu W T_\nu\T\) is linear, so it sends every affine combination of atoms in \(\mathcal K\) to the same affine combination of their images. Since \(\mathcal A_{\nu,e_3}\) is affine and contains each such image, it follows that \(T_\nu(\aff\mathcal K)T_\nu\T\subseteq\mathcal A_{\nu,e_3}\). By Lemma~\ref{lem:circle}, \(\aff\mathcal K=\mathcal B\), and hence \(T_\nu\mathcal B T_\nu\T\subseteq\mathcal A_{\nu,e_3}\). Because \(T_\nu\) is invertible, both affine spaces have dimension \(14\), so the inclusion gives
\[
T_\nu\mathcal B T_\nu\T = \mathcal A_{\nu,e_3}.
\]

For any \(X\in\mathcal A_{\nu,e_3}\) with \(X\succeq0\), set \(W=T_\nu^{-1}X(T_\nu\T)^{-1}\). By the equality just proved we have \(W\in\mathcal B\), and by congruence we have \(W\succeq0\). By Lemma~\ref{lem:circle}, write \(W=\sum_{i=1}^m\lambda_i g_i g_i\T\), where \(g_i g_i\T\in\mathcal K\), \(\lambda_i\ge0\), and \(\sum_{i=1}^m\lambda_i=1\). Substituting this expression into \(X=T_\nu W T_\nu\T\) gives
\[
X = \sum_{i=1}^m\lambda_i\,T_\nu g_i g_i\T T_\nu\T
  = \sum_{i=1}^m\lambda_i\,(T_\nu g_i)(T_\nu g_i)\T.
\]
Each \((T_\nu g_i)(T_\nu g_i)\T\) belongs to \(\Sset\), so this is a convex combination of elements of \(\Sset\) and \(X\in\conv\Sset\).

Thus, for every unit \(w\) and every scalar \(\nu\) with \(|\nu|<1\), Steps~1 and~2 give
\[
\begin{array}{c}
X\in\mathcal A_{\nu,w}\\[-2pt] X\succeq0
\end{array}
\xrightarrow{\ \text{Step 1}\ }
\begin{array}{c}
\widetilde X\in\mathcal A_{\nu,e_3}\\[-2pt] \widetilde X\succeq0
\end{array}
\xrightarrow{\ \text{Step 2}\ }
\widetilde X\in\conv\Sset
\xrightarrow{\ \text{rotate back}\ }
X\in\conv\Sset.
\]
For the last arrow, rotate both columns of each frame in the decomposition by \(P\T\). It remains to find a unit \(w\) and a scalar \(\nu\) with \(|\nu|<1\) satisfying \(L(X)\binom{-\nu}{w}=0\) whenever \(X\in\Sigma\) and \(L(X)\) is singular.

\emph{Step 3: Existence of a latitude kernel vector.}
Let \(X\in\Sigma\) have singular \(L(X)\). By the Schur complement \eqref{eq:LPhi}, \(L(X)\succeq0\) implies \(\Phi(X)\succeq0\). Since \(\tr D=\tr X_{11}+\tr X_{22}=2\), we have
\[
0 \le \tr\Phi(X) = \tr(I_3-D-vv\T) = 1-\norm{v}^2,
\]
so \(\norm{v}\le1\).
\begin{itemize}
\item If \(\norm{v}<1\), choose a nonzero vector \((-\nu;w)\in\ker L(X)\). Block multiplication gives
\[
0 = L(X)\begin{pmatrix}-\nu\\ w\end{pmatrix}
  = \begin{pmatrix}1 & v\T\\ v & I_3-D\end{pmatrix}\begin{pmatrix}-\nu\\ w\end{pmatrix}
  = \begin{pmatrix}-\nu+v\T w\\ -\nu v+(I_3-D)w\end{pmatrix}.
\]
Thus \(\nu=v\T w\), which forces \(w\ne0\) (otherwise \(\nu=0\) as well). Divide the kernel vector \((-\nu;w)\) by \(\norm{w}\), and retain the notation \((-\nu;w)\) for the normalized vector. Then \(\norm{w}=1\), and Cauchy--Schwarz gives \(|\nu|=|v\T w|\le\norm{v}<1\).
\item If \(\norm{v}=1\), then \(\Phi(X)\succeq0\) and \(\tr\Phi(X)=0\) force \(\Phi(X)=0\), i.e., \(I_3-D=vv\T\). Since \(\dim v^\perp=2\), we may choose a unit vector \(w\perp v\), obtaining
\[
L(X)\begin{pmatrix}0\\ w\end{pmatrix}
= \begin{pmatrix}v\T w\\ (I_3-D)w\end{pmatrix}
= \begin{pmatrix}v\T w\\ (v\T w)v\end{pmatrix}
= 0,
\]
so \((0;w)\in\ker L(X)\) with \(\nu=0\in(-1,1)\).
\end{itemize}
In either case \(X\in\mathcal A_{\nu,w}\) for some \(\norm{w}=1\) and \(|\nu|<1\), and Steps~1--2 yield \(X\in\conv\Sset\).
\end{proof}

\subsection{Extreme points and completion of the proof}\label{sec:extreme}

\begin{lemma}[Extreme points have singular \(L\)]\label{lem:extreme}
No matrix \(X\in\Sigma\) with \(L(X)\succ0\) is an extreme point of \(\Sigma\).
\end{lemma}

\begin{proof}
Suppose \(X\in\Sigma\) satisfies \(L(X)\succ0\), equivalently \(\Phi(X)\succ0\) by the Schur complement \eqref{eq:LPhi}. Since \(X\succeq0\) and \(\tr X=2\), we may factor \(X=YY\T\) with \(Y=(Y_1;Y_2)\in\R^{6\times r}\), where \(r=\rank X\ge1\) and \(Y_1,Y_2\in\R^{3\times r}\). Setting \(Z=Y_1+iY_2\in\C^{3\times r}\), we have
\begin{equation}\label{eq:complex-gram}
ZZ^* = (Y_1+iY_2)(Y_1\T-iY_2\T)
= (Y_1Y_1\T+Y_2Y_2\T)+i(Y_2Y_1\T-Y_1Y_2\T)
= D+i[v]_\times.
\end{equation}

Let \(C=I_3-D\). Since \(\Phi(X)=C-vv\T\succ0\), we have \(C\succ vv\T\succeq0\), and \(\tr C=3-\tr D=1\). Thus the eigenvalues \(c_1,c_2,c_3\) of \(C\) are positive and sum to \(1\), which implies \(d_i:=1-c_i>0\) for each \(i\) and hence \(D\succ0\). In an orthonormal eigenbasis of \(C\) (which simultaneously diagonalizes \(D\) and its adjugate \(\adj D\)), for pairwise distinct \(i,j,k\in\{1,2,3\}\) we have
\begin{equation}\label{eq:adj-dom}
(\adj D)_{ii} = d_j d_k = (1-c_j)(1-c_k) 
= 1-c_j-c_k+c_j c_k 
= c_i + c_j c_k > c_i.  
\end{equation}
Therefore \(\adj D\succ C\succ vv\T\).

We claim that \(D+i[v]_\times\succ0\). Take any nonzero \(z=x+iy\in\C^3\) (\(x,y\in\R^3\)) and set \(w=x\times y\). Using skew-symmetry of \([v]_\times\) and \(x\T(v\times y)=-v\T(x\times y)\), we expand
\[
z^*\bigl(D+i[v]_\times\bigr)z
= x\T Dx + y\T Dy - 2x\T[v]_\times y
= x\T Dx + y\T Dy + 2v\T w.
\]
If \(w=0\), then \(z^*(D+i[v]_\times)z=x\T Dx+y\T Dy>0\) since \(D\succ0\) and \(z\ne0\). If \(w\ne0\), then using the scalar triple product, multiplicativity of the determinant, and symmetry of \(D^{-1/2}\) in the last step, we have for every \(u\in\R^3\),
\begin{align*}
u\T\bigl((D^{1/2}x)\times(D^{1/2}y)\bigr)
&= \det[u,D^{1/2}x,D^{1/2}y] 
= (\det D^{1/2})\det[D^{-1/2}u,x,y]\\
&= (\det D^{1/2})(D^{-1/2}u)\T w
 = u\T\bigl((\det D^{1/2})D^{-1/2}w\bigr),
\end{align*}
so \((D^{1/2}x)\times(D^{1/2}y)=(\det D^{1/2})D^{-1/2}w\). Taking squared norms and applying Lagrange's identity yields
\[
w\T(\adj D)w = (\det D)\,w\T D^{-1}w
= \norm{(D^{1/2}x)\times(D^{1/2}y)}^2
= (x\T Dx)(y\T Dy) - (x\T Dy)^2.
\]
Combining \(\adj D\succ vv\T\) (which gives strict inequality since \(w\ne0\)) with this identity and the AM--GM inequality gives
\begin{equation}\label{eq:quad-pos}
2|v\T w| < 2\sqrt{w\T(\adj D)w}  \le 2\sqrt{(x\T Dx)(y\T Dy)}  
\le x\T Dx + y\T Dy,  
\end{equation}
so \(z^*(D+i[v]_\times)z \ge x\T Dx+y\T Dy-2|v\T w|>0\). Thus \(z^*(D+i[v]_\times)z>0\) in both cases.

By \eqref{eq:complex-gram}, \(ZZ^*=D+i[v]_\times\) is a positive definite \(3\times3\) matrix and hence has full rank. Since \(\im(ZZ^*)\subseteq\im Z\) and \(Z\) has three rows and \(r=\rank X\) columns, we obtain
\[
3 = \rank(ZZ^*) \le \rank Z \le \min\{3,r\},
\]
which forces \(r\ge3\).

Finally, suppose for contradiction that \(X\) is an extreme point of \(\Sigma\). The set \(U=\{Y\in\mathcal A:L(Y)\succ0\}\) is open in \(\mathcal A\) and contains \(X\). Set \(\mathcal F=\{Y\in\mathcal A:Y\succeq0\}\). Suppose \(X\) is not an extreme point of \(\mathcal F\). Then \(X\) lies strictly between two distinct points of \(\mathcal F\), so there is an \(E\ne0\) with \(X\pm E\in\mathcal F\). For sufficiently small \(t\in(0,1)\), both \(X\pm tE\) lie in \(\mathcal F\) by convexity and in \(U\) by openness. They are distinct elements of \(\mathcal F\cap U\subseteq\Sigma\), and
\[
X = \tfrac12(X+tE) + \tfrac12(X-tE).
\]
This contradicts the assumption that \(X\) is an extreme point of \(\Sigma\). Thus \(X\) is an extreme point of \(\mathcal F\). Because \(\mathcal F\) is defined by the three linear equalities \eqref{eq:affineA}, the Barvinok--Pataki rank bound \cite[Theorem~2.1]{Pataki} (see also \cite{Barvinok}) gives \(r(r+1)/2\le3\), i.e., \(r\le2\), contradicting \(r\ge3\).
\end{proof}

\begin{proof}[Proof of Theorem~\ref{thm:main}]
For any lifted frame \(qq\T\in\Sset\) with \(q=(a;b)\), \(Q=[a,b]\in\St(3,2)\), and \(n=a\times b\), we have \(qq\T\succeq0\); the orthonormality \(Q\T Q=I_2\) gives the trace constraints \eqref{eq:affineA}; and \(aa\T+bb\T+nn\T=I_3\) together with \eqref{eq:cross} gives \(L(qq\T)=(1;n)(1;n)\T\succeq0\). Hence \(\Sset\subseteq\Sigma\), and convexity of \(\Sigma\) yields \(\conv\Sset\subseteq\Sigma\) (see also \cite{DSX}).

For the reverse inclusion, first note that \(\Sigma\) is compact: it is closed because \(\mathcal A\) and the positive semidefinite cones are closed and \(L\) is continuous, and it is bounded because \(X\succeq0\) and \(\tr X=\tr X_{11}+\tr X_{22}=2\) imply \(\norm{X}_{\mathrm F}\le\tr X=2\). By Minkowski's theorem \cite[Corollary~18.5.1]{Rockafellar}, \(\Sigma\) is the convex hull of its extreme points. Every extreme point of \(\Sigma\) has singular \(L(X)\) by Lemma~\ref{lem:extreme}, and therefore lies in \(\conv\Sset\) by Lemma~\ref{lem:latitude}. Thus \(\Sigma\subseteq\conv\Sset\).
\end{proof}

\section{Conclusion}

We have established the exact semidefinite description of the
quadratic convex hull of real \(3\times2\) Stiefel matrices, resolving
the conjecture of Dodig, Sto\v{s}i\'c, and Xavier. The result shows
that the lifted orthogonality constraints, together with the
normal-moment constraint, suffice to represent every feasible
matrix as a convex combination of lifted frames. It follows that
the proposed semidefinite formulation is exact for every
homogeneous quadratic objective over \(\St(3,2)\).
A natural direction for future work is to develop an explicit
procedure for recovering a finite frame decomposition from a
feasible matrix of the semidefinite formulation.

\medskip
\noindent\textbf{Funding.} \\
This work was supported by the National Natural Science Foundation of China (Grant No.~12631012).

\medskip
\noindent\textbf{Disclosure.}\\ 
During the preparation of this manuscript, the authors used
AI-assisted tools for language editing, formatting, and assistance
in checking mathematical arguments. The authors
independently verified all mathematical statements and proofs
and take full responsibility for the final manuscript. 

{\small\emergencystretch=2.5em\setlength{\itemsep}{3pt plus 1pt}

}
\end{document}